\documentclass[reqno]{amsart}
 \usepackage{amsbsy,amssymb,amscd,amsfonts,latexsym,amstext,delarray,
 amsmath,graphicx, color}
 \usepackage[dvips]{epsfig}
\usepackage{hyperref}
\input xypic

\newtheorem{thm}{Theorem}[section]
\newtheorem{prop}[thm]{Proposition}

\newtheorem{lem}[thm]{Lemma}

\newtheorem{defn}[thm]{Definition}

\def\R{{\mathbb R}}
\def\Z{{\mathbb Z}}

\def\T{{\mathbb T}}

\def\cD{{\mathcal D}}

\def\cH{{\mathcal H}}

\def\cL{{\mathcal L}}

\def\tr{{\rm Trace}}

\newcommand{\ie}{{\it i.e.\/}\ }

\newcommand{\cf}{{\it cf.\/}\ }

\def\qqq{\,,\quad~\forall}

\def\ries{\Sigma}
\def\lap{\triangle}
\def\modu{\Delta}
\def\pert{\triangle'}

\newcommand{\nil}[1]{}

\title{The Gauss-Bonnet Theorem for  the noncommutative two torus}
\author{Alain CONNES and Paula TRETKOFF}

\begin{document}

\maketitle

\begin{abstract}
In this paper we shall show that the value at the origin, $\zeta
(0)$, of the zeta function of the Laplacian on the non-commutative
two torus, endowed with its canonical conformal structure, is
independent of the choice of the volume element (Weyl factor) given
by a (non-unimodular) state. We had obtained, in the late eighties,
in an unpublished computation,  a general formula for $\zeta (0)$
involving modified logarithms of the modular operator of the state.
We give here the detailed computation and prove that the result is
independent of the Weyl factor as in the classical case, thus
proving the analogue of the Gauss-Bonnet theorem for  the
noncommutative two torus.
\end{abstract}

 \tableofcontents

\section{Introduction}\label{sec1}

The main result of this paper, namely the analogue of the
Gauss-Bonnet theorem for  the noncommutative two torus
$\T^2_\theta$, follows from the computation of the  value at the
origin, $\zeta (0)$, of the zeta function  of the analogue of the
Laplacian for non-unimodular geometric structures on $\T^2_\theta$.
The two authors did this computation at the end of the $1980$'s. The
 result was mentioned in \cite{paula}, but it was never published except
as an MPI preprint \cite{cc}. At the time the computation was
done, the result was formulated in terms of the modular operator
but its significance was unclear. In the mean time the following two
theories have been developed:
\begin{itemize}
  \item The spectral action
  \item The non-unimodular spectral triples
\end{itemize}
The spectral action \cite{cc2}, \cite{mc2}, allows one to interpret
gravity coupled to the standard model from spectral invariants of a
geometry which encodes a fine structure of space-time. In the two
dimensional case the constant term in the spectral action is (up to
adding the dimension of the kernel of the operator) given by the
value at the origin, $\zeta (0)$, of the zeta function. This value
is a topological invariant in classical Riemannian geometry.

The twisted (or non-unimodular) spectral triples appeared very
naturally \cite{cm} in the study of type III examples of foliation
algebras.

Thus, in more modern terminology, our computation was that of the
spectral action for  spectral triples on the non-commutative two
torus $\T^2_\theta$ both in the usual and the twisted cases.
Initially, the complexity of the computation and the lack of simplicity of the
result made us reticent about publishing it. It is only recently that,
by pushing the computation further, we could prove the expected (\cf
\cite{cc3}) conformal invariance.

\begin{thm}
-- Let $\theta$ be an irrational number and consider, on the
non-commutative two torus $\T^2_\theta$,
 a translation invariant conformal structure. Let $k$ be an invertible positive element of $C^\infty(\T^2_\theta)$ considered as a Weyl factor rescaling the metric.
Then the value at the origin of the zeta function $\zeta (s)$ of the
Laplacian of the rescaled metric is independent of $k$.
\end{thm}

In fact we must retract a bit after stating this theorem since we
 have only performed the computation for the simplest translation
invariant conformal structure but we do not expect that the general case
will be different.

There are two main actors in the computation and they display an
interesting interplay between two theories which look a priori quite
distinct, but use the same notation for their key ingredient:
\begin{itemize}
  \item The spectral theory of the Laplacian $\lap$
  \item The modular operator $\modu$  of states on operator algebras
\end{itemize}
The second is part of the scene because of the non-unimodularity of
the  spectral triple, or in simpler terms because the state defining
the volume form is no longer assumed to be  a trace and hence
inherits a modular operator $\modu$. While both the Laplacian and
the modular operator will be denoted by the capital letter delta, we
hope the distinction\footnote{We use a slightly different typography
for the two cases} between these two operators will remain clear
throughout to the alert reader.

\section{Preliminaries}\label{prel}

Recall that for a classical Riemann surface $\ries$ with metric $g$,
to the Laplacian $\lap_g = d^* d$, where $d$ is the de-Rham
differential operator acting on functions on the Riemann surface, one associates
the zeta function
\begin{equation*}
\zeta (s) = \sum \, \lambda_j^{-s} \, , \quad {\rm Re} (s) > 1 \, ,
\end{equation*}
where the summation is over the non-zero eigenvalues $\lambda_j$ of
$\lap_g$ (see, for example, \cite{Rosen}). The meromorphic continuation of $\zeta (s$) to $s=0$,
where it has no pole, gives the important information
\begin{equation*}
\zeta (0)+{\rm Card} \{ j \mid \lambda_j = 0 \} = \frac{1}{12 \pi}
\int_\ries \, R  = \frac 16 \,\chi (\ries)  \, ,
\end{equation*}
where $R$ is the scalar curvature and $\chi (\ries) $ the
 Euler-Poincar\'e characteristic. This vanishes when $\ries$ is the
classical 2-torus ${\mathbb R}^2 / {\mathbb Z}^2$, for example, and
is an invariant within the conformal class of the metric, that is
under the transformation $g \to e^f g$ for $f$ a smooth real valued
function on $\ries$.

\subsection{The noncommutative two torus $\T^2_\theta$} \hfill\medskip

Fix a real irrational number $\theta$. We consider the irrational
rotation $C^*$-algebra $A_{\theta}$, with two unitary generators
which satisfy
\begin{equation*}
VU = e^{2\pi i \theta} \, UV \, , \quad U^* = U^{-1} \, , \ V^* =
V^{-1} \, .
\end{equation*}
 We introduce the dynamical system given by the action of $\T^2$, $\T = \R/2\pi\Z$ on $A_{\theta}$ by the 2-parameter group of automorphisms $\{ \alpha_s \}$, $s\in \R^2$ determined by
\begin{equation*}
\alpha_s(U^n \, V^m)=e^{is.(n,m)}U^n \, V^m \, , \quad (s \in \R^2)
\, ,
\end{equation*}
We define the sub-algebra $A_{\theta}^{\infty}$ of smooth elements
of $A_{\theta}$ to be those $x$ in $A_{\theta}$ such that the
mapping
\begin{equation*}
{\mathbb R}^2 \to A_{\theta}, \  \ s \mapsto \alpha_s  (x)
\end{equation*}
is smooth. Expressed as a condition on the coefficients of the
element $a\in A_{\theta}$,
 $$
 a=\sum a (n,m)U^nV^m
 $$
 this is the same as saying that they be of rapid decay, namely that $\{ \vert n \vert^k \, \vert m \vert^q \, \vert a (n,m) \vert \}$ be bounded for any positive $k, q$. The derivations associated to the above group of automorphisms are given
by  the defining relations,
\begin{equation*}
\delta_1 (U) =   \, U \, , \quad \delta_1 (V) = 0 \, ,
\end{equation*}
\begin{equation*}
\delta_2 (U) = 0 \, , \quad \delta_2 (V) =   \, V \, .
\end{equation*}
The derivations $\delta_1$, $\delta_2$ are analogues of the
differential operators $ \frac 1i\partial / \partial x$, $ \frac
1i\partial / \partial y$ on the smooth functions on ${\mathbb R}^2 /
{2\pi\mathbb Z}^2$.

\subsection{Conformal structure on $\T^2_\theta$} \hfill\medskip

As $\theta$ is supposed irrational, there is a unique trace $\tau$
on $A_{\theta}$ determined by the orthogonality properties
\begin{equation}\label{trace}
\tau (U^n \, V^m) = 0 \quad \mbox{if} \quad (n,m) \ne (0,0) \, ,
\quad \mbox{and} \quad \tau (1) = 1 \, .
\end{equation}
We can construct a Hilbert space ${\cH_0}$ from $A_{\theta}$ by
completing with respect to the inner product
\begin{equation}\label{inprod}
\langle a,b \rangle = \tau (b^* a) \, , \quad a,b \in A_{\theta} \,
,
\end{equation}
and, using the derivations $\delta_1 , \delta_2$, we introduce a
complex structure by defining
\begin{equation}\label{deltas}
\partial = \delta_1 + i \delta_2 \, , \quad \partial^* =\delta_1 - i \delta_2
\end{equation}
where (extending $\partial , \partial^*$ to unbounded operators on
${\cH_0}$) $\partial^*$ is the adjoint of $\partial$ with respect to
the inner product defined by $\tau$. As an appropriate analogue of
the space of $(1,0)$-forms on the classical 2-torus,  one takes the
unitary bi-module $\cH^{(1,0)}$ over $A_{\theta}^{\infty}$ given by
the Hilbert space completion of the space of finite sums $\sum \, a
\ \partial b$, $a,b \in A_{\theta}^{\infty}$, with respect to the
inner product
\begin{equation}\label{onezero}
\langle a \partial b , a' \partial b' \rangle = \tau ((a')^* \,
a(\partial b) (\partial b')^*) \, , \quad a,a',b,b' \in
A_{\theta}^{\infty} \, .
\end{equation}
The information on the conformal structure is encoded by the
positive Hochschild two cocycle (\cf \cite{concun} \cite{Co-book})
on $A_{\theta}^{\infty}$ given by
\begin{equation}\label{poshoch}
   \psi(a,b,c)=-\tau(a\partial b \partial^* c)
\end{equation}

\subsection{Modular automorphism $\modu$} \hfill\medskip

In order to vary inside the conformal class of a metric we consider
the family of positive linear functionals  $\varphi = \varphi_h$,
parameterized by $h=h^* \in A_{\theta}^{\infty}$, a  self-adjoint
element of $A_{\theta}^{\infty}$, and defined on $A_{\theta}$ by
\begin{equation*}
\varphi (a) = \tau (ae^{-h}) \, , \quad a \in A_{\theta} \, .
\end{equation*}
Note that, whereas for $\tau$ we have the trace relation
\begin{equation*}
\tau (b^* a) = \tau (ab^*) \, , \quad a,b \in A_{\theta} \, ,
\end{equation*}
for $\varphi$ we have
\begin{equation*}
\varphi (ab) = \varphi (b e^{-h} ae^{h}) = \varphi (b \sigma_i \,
(a)) \, , \quad a \in A_{\theta} \, ,
\end{equation*}
which is the KMS condition at $\beta=1$ for the 1-parameter group
$\sigma_t$, $t \in {\mathbb R}$, of inner automorphisms
\begin{equation*}
\sigma_t (x) = e^{ith} x e^{-ith}.
\end{equation*}
This group is $\sigma_t=\modu^{-it}$  where the modular operator is
\begin{equation*}
\modu (x) = e^{-h} x e^h
\end{equation*}
which is a positive operator fulfilling
\begin{equation}\label{modoper}
    \langle \modu^{1/2}x, \modu^{1/2}x  \rangle_\varphi=\langle x^*,x^*   \rangle_\varphi\qqq x\in A_{\theta}.
\end{equation}
where we define the inner product $ \langle \ , \ \rangle_\varphi$
on $A_{\theta}$ by
\begin{equation*}
\langle a,b\rangle_\varphi = \varphi (b^* a) \, , \quad a,b \in
A_{\theta} \, .
\end{equation*}
We let $\cH_\varphi$ be the Hilbert space completion of $A_{\theta}$
for the inner product $ \langle \ , \ \rangle_\varphi$. It is a
unitary left module on $A_{\theta}$ by construction. The 1-parameter
group $\sigma_t$ is generated by the derivation  $- \log \modu$
\begin{equation*}
-\log \modu(x) = [h,x] \, , \quad x \in A_{\theta}^{\infty} \, .
\end{equation*}

\subsection{Laplacian and Weyl factor on $\T^2_\theta$} \hfill\medskip

The Laplacian $\lap$ on functions on $\T^2_\theta$ is given by
\begin{equation}\label{lap0}
\lap = \partial^* \partial = \delta_1^2 + \delta_2^2 \, ,
\end{equation}
where $\partial$ is viewed as an unbounded operator from $\cH_0$ to
$\cH^{(1,0)}$.

When one modifies the volume form on $\T^2_\theta$ by replacing the
trace $\tau$ by the functional $\varphi$ the modified Laplacian
$\pert$ is given by
\begin{equation}\label{lap1}
\pert = \partial_\varphi^* \partial_\varphi
\end{equation}
where $\partial_\varphi$ is the same operator $\partial$ but viewed
as an unbounded operator from $\cH_\varphi$ to $\cH^{(1,0)}$. By
construction $\pert$ is a positive unbounded operator in
$\cH_\varphi$.

\begin{lem}\label{laplem}
The operator $\pert$ is anti-unitarily equivalent to the positive
unbounded operator $k\lap k$ in the Hilbert space $\cH_0$, where
$k=e^{h/2}\in A_{\theta}$ acts in $\cH_0$ by left multiplication.
\end{lem}
\proof The right multiplication by $k$ extends to an isometry $W$
from $\cH_0$ to $\cH_\varphi$,
\begin{equation}\label{lap2}
    W a=ak\qqq a\in A_{\theta}
\end{equation}
since
$$
\langle Wa,Wb\rangle_\varphi=\tau((bk)^*(ak)k^{-2})=\tau(b^*a)\qqq
a,b \in A_{\theta}.
$$
The operator $ \partial_\varphi\circ W$ from $\cH_0$ to
$\cH^{(1,0)}$ is given by
$$
\partial_\varphi\circ W(a)=\partial (ak)=\partial\circ R_k
$$
where $R_k$ is the right multiplication by $k$. Thus $\pert =
\partial_\varphi^* \partial_\varphi$ is unitarily equivalent to
$(\partial\circ R_k)^*\partial\circ R_k=R_k\partial^* \partial R_k$.
Now, let $J$ be the anti-unitary involution on $\cH_0$ given by the
star operation $Ja=a^*$ for all $a\in A_{\theta}$. The operator $J$
commutes with $\lap$ and fulfills $JR_kJ=k$ using $k^*=k$. This gives
 the required equivalence. \endproof

It is the dependence on $k$ in the computations of the behavior of
the zeta function of $k\lap k$ at the origin that will feature in
what follows.

\subsection{Spectral Triples} \hfill\medskip

With the notations of the previous section, consider the Hilbert
space and operator
\begin{equation}\label{spectrip}
\cH=\cH_\varphi\oplus \cH^{(1,0)}\,, \ \ D= \left(
  \begin{array}{cc}
    0 & \partial^*_\varphi \\
    \partial_\varphi & 0\\
  \end{array}
\right).
\end{equation}
We recall that $\cH^{(1,0)}$ is naturally a unitary bimodule over
$A_{\theta}$.

\begin{lem}\label{speclem}
1) The left action of $A_{\theta}$ in $\cH=\cH_\varphi\oplus
\cH^{(1,0)}$ and the operator $D$ yield an even spectral triple
$(A_{\theta},\cH,D)$.

2) Let $J_\varphi$ be the Tomita antilinear unitary in $\cH_\varphi$
and $a\mapsto J_\varphi a^* J_\varphi$ the corresponding unitary
right action of $A_{\theta}$ in $\cH_\varphi$. Then the right action
$a\mapsto a^{\rm op}$ of $A_{\theta}$ in $\cH=\cH_\varphi\oplus
\cH^{(1,0)}$ and the operator $D$ yield an even twisted spectral
triple $(A^{\rm op}_{\theta},\cH,D)$, \ie the following operators
are bounded
\begin{equation}\label{rightbounded}
    D \,a^{\rm op}- (k^{-1}ak)^{\rm op} D \qqq a\in A_{\theta}.
\end{equation}
3) The zeta function of the operator $D$ \ie
$\zeta_D(s)=\tr(|D|^{-s})$ is equal to $2\tr((k\lap k)^{-s/2})$.
\end{lem}

\proof 1) In order to show that $[D,a]$ is bounded, it is enough to
check that $[\partial_\varphi,a]$ is bounded which follows from the
derivation property of $\partial_\varphi$ and the equivalence of the
norms $\Vert.\Vert_\varphi$ and $\Vert.\Vert_0$.

2) Let us first show that the twisted commutator $\partial_\varphi
\,a^{\rm op}-(k^{-1}ak)^{\rm op}\,\partial_\varphi$ is bounded or
equivalently that
\begin{equation}\label{bound1}
\partial_\varphi \,(kak^{-1})^{\rm op}-a^{\rm op}\,\partial_\varphi
\end{equation}
is bounded as an operator from $\cH_\varphi$ to $\cH^{(1,0)}$. Using the isometry
$W$ from $\cH_0$ to $\cH_\varphi$ defined in \eqref{lap2}, replaces
$\partial_\varphi$ by $\partial\circ R_k$ and replaces $a^{\rm op}$
by the right multiplication $R_a$ by $a$ for any $a\in A_{\theta}$.
Thus it replaces the operator \eqref{bound1} by
$$
\partial\circ R_k\circ R_{kak^{-1}}-R_a\circ \partial\circ R_k=R_{\partial a}\circ R_k
$$
which is a bounded operator from $\cH_0$ to $\cH^{(1,0)}$. Thus
\eqref{bound1} is bounded and so is its adjoint
$$
((kak^{-1})^{\rm op})^*\partial_\varphi^*-\partial_\varphi^*(a^{\rm
op})^*.
$$
Thus the boundedness of \eqref{rightbounded} follows from the
equality
$$
((kak^{-1})^{\rm op})^*=((kak^{-1})^*)^{\rm op}=(k^{-1}a^*k)^{\rm
op}
$$
so that $\partial_\varphi^* \,a^{\rm op}-(k^{-1}ak)^{\rm
op}\,\partial_\varphi^*$ is bounded for all $a\in A_{\theta}$.

3) We have seen in Lemma \ref{laplem} that the spectrum of $\pert =
\partial_\varphi^* \partial_\varphi$ is the same as that of $k\lap
k$. Since the non-zero part of the spectrum of a product $A^*A$ is
the same as the non-zero part of the spectrum of $AA^*$, using the
unitary equivalence given by the polar decomposition, one gets the
required result. \endproof

\bigskip

\section{Statement of the theorem}\label{sec2}

With the notation of \S\ref{sec1}, we study the Laplacian zeta
function defined for ${\rm Re} (s) > 1$ by the Mellin transform
\begin{equation*}
\zeta (s) = \frac{1}{\Gamma (s)} \int_0^{\infty} {\rm Trace}^+
(e^{-t\pert}) \, t^{s-1} \, dt = {\rm Trace} (\pert^{-s}) \, .
\end{equation*}
Here, $\pert\sim k \, \lap \, k$, and
\begin{equation*}
{\rm Trace}^+ (e^{-t\pert}) = {\rm Trace} \, (e^{-t\pert}) - {\rm
Dim} \, {\rm Ker} (\pert) \, ,
\end{equation*}
where by ${\rm Trace} (\cdot)$ we understand the ordinary trace of
the operator. The definition of $\zeta (s)$ can be extended by
meromorphic continuation to all values of $s$, barring $s = 1$ where
the function has a simple pole. We now state the main result.

\begin{thm}\label{thmmain}
-- Let $\theta$ be an irrational number and $k$ an invertible
positive element of $A_{\theta}^{\infty}$. Then the value at the
origin of the zeta function $\zeta (s)$ of the operator $\pert \sim
k \, \lap \, k$ is independent of $k$.
\end{thm}

\bigskip

Our proof relies on a long explicit computation whose main
ingredient is the following lemma:

\subsection{Main technical Lemma} \hfill\medskip

\begin{lem}\label{main}
-- Let $\theta$ be an irrational number and $k$ an invertible
positive element of $A_{\theta}^{\infty}$. Then the value at the
origin of the zeta function $\zeta (s)$ of the operator $\pert\sim k
\, \lap \, k$ is given by\footnote{summation of repeated indices
over $j = 1,2$ is understood}
\begin{equation}\label{main1}
\zeta (0)+1 =2\pi\, \varphi (f(\modu)(\delta_j (k))\delta_j (k))
\end{equation}
where $\varphi$ is the functional $\varphi(x)=\tau(xk^{-2})$,
$\modu$ is the modular operator of $\varphi$ and the function $f(u)$
is given by
\begin{equation}\label{fu}
   f(u)=\frac 16 u^{-1/2}-\frac 13+{\mathcal L}_1 (u)-2(1 + u^{1/2}){\mathcal L}_2 (u) +(1 + u^{1/2})^2 \, {\mathcal L}_3 (u)
\end{equation}
where ${\mathcal L}_m$, $m$ a positive integer, stands for the
modified logarithm
\begin{equation*}
{\mathcal L}_m (u)= (-1)^m \, (u - 1)^{-(m+1)} \left( \log u -
\sum_{j=1}^m (-1)^{j+1} \frac{ (u - 1)^j}{j} \right) \, .
\end{equation*}
\end{lem}

The proof of this lemma will occupy the remaining sections of the
paper. In the present section, we shall show how it implies Theorem
\ref{thmmain}. The statement of Lemma \ref{main} is the
same\footnote{up to adding $1$ to $\zeta(0)$} as that of the main
Theorem in \cite{cc}, where the right hand side is written as
$\tau(h(\theta,k))$ where
\begin{equation}\label{former}
h(\theta,k)=\frac \pi 3 k^{-1}\delta_j^2 (k)-\frac{2\pi}{3}k^{-1}
\delta_j (k)\delta_j (k)    k^{-1}+2\pi \cD_1(k^{-1} \delta_j
(k))\delta_j (k)    k^{-1}
\end{equation}
 $$-4\pi \cD_2(1+\modu^{1/2})(k^{-1} \delta_j (k))\delta_j (k)    k^{-1}+2\pi \cD_3(1+2\modu^{1/2}+\modu)(k^{-1} \delta_j (k))\delta_j (k)    k^{-1}
$$
where $\cD_n=\cL_n(\modu)$. One has indeed
$$
\tau(k^{-1}\delta_j^2 (k))=-\tau(\delta_j(k^{-1})\delta_j
(k))=\tau(k^{-1}\delta_j (k)k^{-1}\delta_j (k))
=\tau(k^{-2}\modu^{-1/2}(\delta_j (k))\delta_j (k))
$$
which accounts for the first term in the right hand side of
\eqref{fu}, while the other terms are easily compared since the left
multiplication by $k^{-1}$ commutes with any function of $\modu$.

\subsection{Result in terms of $\log(k)$} \hfill\medskip

The function $f(u)$ is of the form $h(\log(u))$ where $h$ is the
entire function
\begin{equation}\label{main2}
h(x)=-\frac{e^{-x/2} \left(-1+3 e^{x/2}+3 e^x+6 e^{3 x/2} x-3 e^{2
x}-3 e^{5 x/2}+e^{3 x}\right)}{6 \left(-1+e^{x/2}\right)^4
\left(1+e^{x/2}\right)^2}.
\end{equation}
Thus one can rewrite \eqref{main1} in the form
\begin{equation}\label{main3}
\zeta (0) +1 =2\pi\, \varphi (h(\log\modu)(\delta_j (k))\delta_j
(k))
\end{equation}

As pointed out in \S\ref{sec1}, in the commutative case the
corresponding value for the zeta function at the origin is zero. In
fact in that case one has $\log\modu=0$. The function $h$ vanishes
at $0$ and its Taylor expansion there is
$$
h(x)=-\frac{x}{20}+\frac{x^2}{40}-\frac{x^3}{210}+\frac{x^4}{3360}+\frac{x^5}{201600}+O[x]^6.
$$

 \begin{figure}
\begin{center}
\includegraphics[scale=0.9]{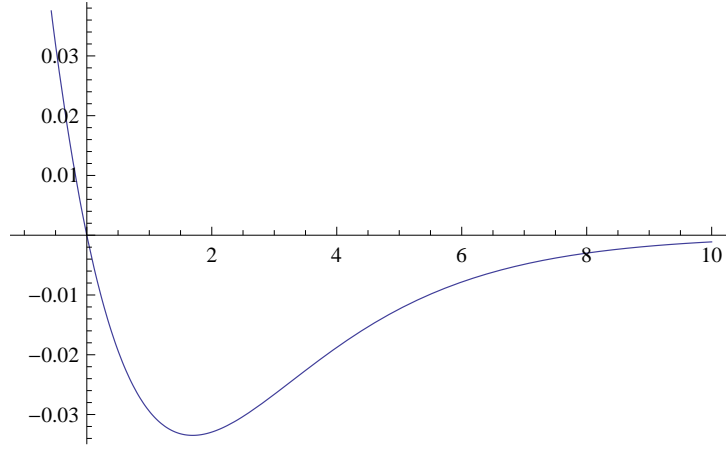}
\end{center}
\caption{Graph of the function $h$.\label{funchdraw} }
\end{figure}

In fact, one obtains a further simplification in general by
expressing the result in terms of the element $\psi=\log k$. We
introduce the function
\begin{equation}\label{fctK}
K(x)=-\frac{x- \text{sh}\left[\frac{x}{2}\right]- \text{sh}[x]+\frac
13\text{sh}\left[\frac{3 x}{2}\right]}{
x^2\text{sh}\left[\frac{x}{2}\right]^2}.
\end{equation}

\begin{lem}\label{cormain}
-- With the notations of Theorem \ref{main}, one has
\begin{equation}\label{cormain1}
\zeta (0) +1 =2\pi\, \tau (K(\log\modu)(\delta_j (\log k))\delta_j
(\log k)).
\end{equation}
\end{lem}

\proof One uses the following formula
\begin{equation}\label{deltajk}
    k^{-1}\delta_j(k)=\int_0^1\modu^{s/2}(\delta_j (\log k))ds
\end{equation}
which gives
\begin{equation}\label{deltajkbis}
    k^{-1}\delta_j(k)=2\,\frac{\modu^{1/2}-1}{\log \modu}(\delta_j (\log k))
\end{equation}
and similarly
$$
\delta_j(k)k^{-1}=-2\,\frac{\modu^{-1/2}-1}{\log \modu}(\delta_j
(\log k)).
$$
One then rewrites \eqref{main1} in the form
$$
\zeta (0) +1 =2\pi\, \tau (f(\modu)(k^{-1}\delta_j (k))\delta_j
(k)k^{-1})
$$
and one uses the following identity, where $F$ is an entire function,
\begin{equation}\label{byparts}
    \tau (a F(\log\modu)(b))=\tau (F(-\log\modu)(a)b)\qqq a,b \in A_{\theta}^{\infty}
\end{equation}
which shows that \eqref{cormain1} holds with
$$
K(x)=4\left(-1+e^{x/2}\right)^2 \,x^{-2}\,h(x)
$$
which simplifies to \eqref{fctK}. \endproof

\subsection{Proof of Theorem \ref{thmmain}} \hfill\medskip

By \eqref{fctK} the function $K(x)$ is an odd function. Thus using
\eqref{byparts} one has
$$
\tau (K(\log\modu)(\delta_j (\log k))\delta_j (\log k))=-\tau
(K(\log\modu)(\delta_j (\log k))\delta_j (\log k))=0.
$$

\section{Pseudo-differential calculus}\label{sec3}

With the notation of the preceding sections, we introduce in the
present one the notion of a pseudo-differential operator given the
dynamical system  $(A_{\theta}^{\infty} , \alpha_s)$
 as developed in \cite{B} and \cite{C}. First of all, for a non-negative integer $n$, we define the vector space of differential operators of order at most $n$ to be those polynomial expressions in $\delta_1 , \delta_2$ of the form
\begin{equation*}
P(\delta_1 , \delta_2) = \sum_{\vert j \vert \leq n} a_j \,
\delta_1^{j_1} \, \delta_2^{j_2} \, , \quad a_j \in
A_{\theta}^{\infty} \, , \quad j = (j_1 , j_2) \in {\mathbb Z}_{\geq
0}^2 \, , \quad \vert j \vert = j_1 + j_2 \, .
\end{equation*}
To extend this definition, let $\R^2$ be the group dual to ${\mathbb
R}^2$ and introduce the class of operator valued distributions given
by those complex linear functions $P : C^{\infty} (\R^2) \to
A_{\theta}^{\infty}$ which are continuous with respect to the
semi-norms $p_{i_1 , i_2}$ determined by
\begin{equation*}
p_{i_1 , i_2} (P(\varphi)) = \Vert \delta_1^{i_1} \, \delta_2^{i_2}
(P(\varphi)) \Vert \, , \quad i_1 , i_2 \in {\mathbb Z}_{\geq 0} \,
, \quad \varphi \in C^{\infty} (\R^2) \, .
\end{equation*}
We use the notation $y_1 = e^{2\pi i \xi_1}$, $y_2 = e^{2\pi i
\xi_2}$, $\xi = (\xi_1 , \xi_2) \in {\mathbb R}^2$, for the
canonical coordinates of $\R^2$, and $\partial_1 = \partial /
\partial \xi_1$, $\partial_2 = \partial / \partial \xi_2$ for the
corresponding derivations. We may now introduce the algebra of
pseudo-differential operators via the algebra of operator valued
symbols.

\bigskip

\begin{defn}
-- An element $\rho = \rho (\xi) = \rho (\xi_1 , \xi_2)$ of
$C^{\infty} (\R^2 , A_{\theta}^{\infty})$ is a symbol of order the
integer $n$ if and only if for all non-negative integers $i_1 , i_2
, j_1 , j_2$
\begin{equation*}
p_{i_1 , i_2} (\partial_1^{j_1} \, \partial_2^{j_2} \, \rho (\xi))
\leq c \, (1+\vert \xi \vert)^{n-\vert j \vert} \, ,
\end{equation*}
where $c$ is a constant depending only on $\rho$, and if there
exists an element $k = k (\xi_1 , \xi_2)$ of $C^{\infty} (\R^2 - \{
0,0 \} , A_{\theta}^{\infty})$ such that
\begin{equation*}
\lim_{\lambda \to \infty} \lambda^{-n} \, \rho (\lambda \, \xi_1 ,
\lambda \, \xi_2) = k(\xi_1 , \xi_2) \, .
\end{equation*}
\end{defn}

\bigskip

We denote the space of symbols of order $n$ by $S_n$, the union $S =
\cup_{n \in {\mathbb Z}} \, S_n$ forming an algebra. Symbols of
non-integral order are not required for this paper. An example of a
symbol of order $n$ a positive integer is provided by the polynomial
$\rho (\xi) = \sum_{\vert j \vert \leq n} \, a_j  \, \xi_1^{j_1} \,
\xi_2^{j_2}$, $a_j \in A_{\theta}^{\infty}$, and one has $\rho(n,m)
= \sum_{\vert j \vert \leq n} \, a_j \, n^{j_1} \, m^{j_2}$ so that
$\rho(n,m) \, U^n \, V^m = \sum_{\vert j \vert \leq n} \, a_j \,
\delta_1^{j_1} \, \delta_2^{j_2} (U^n \, V^m)$. For an element $a =
\sum_{n,m}$ $a(n,m) \, U^n \, V^m$ of $A_{\theta}^{\infty}$ one
therefore has $\sum_{n,m} \, \rho(n,m) \, a(n,m) \, U^n \, V^n =
\sum_{\vert j \vert \leq n}$ $a_j \, \delta_1^{j_1} \,
\delta_2^{j_2} (a)$, associating to the symbol $\rho$ the
differential operator $P_{\rho} = P(\delta_1 , \delta_2)$ $=
\sum_{\vert j \vert \leq n} \, a_j \, \delta_1^{j_1} \,
\delta_2^{j_2}$ on $A_{\theta}^{\infty}$.

 For every integer $n$, a symbol $\rho$ of that order determines an operator on $A_{\theta}^{\infty}$ via the map $\psi : \rho \mapsto P_{\rho}$ given by the general formula
 \begin{equation}\label{defpsido}
P_{\rho} (a) = (2\pi)^{-2}\int\int\,e^{-is.\xi} \rho(\xi) \,
\alpha_s(a)dsd\xi \, .
\end{equation}
In our case this gives, using
$$
\alpha_s(U^n \, V^m)=e^{is.(n,m)}U^n \, V^m
$$
the simpler formula
\begin{equation}\label{defpsido1}
P_{\rho} (a) = \sum_{n,m \in {\mathbb Z}} \, \rho(n,m) \, a(n,m) \,
U^n \, V^m \, , \quad a = \sum_{n,m} \, a(n,m) \, U^n \, V^m \, .
\end{equation}
For example, the image under $\psi$ of the symbol $(1+\vert \xi
\vert^2)^{-k}$, $k \geq 1$, of order $-2k$ acts on
$A_{\theta}^{\infty}$.

\bigskip

\begin{defn}
-- The space $\psi$ of pseudo-differential operators is given by the
image of the algebra $S$ under the map $\psi$.
\end{defn}

\bigskip

\begin{defn}
-- The equivalence $\rho \sim \rho'$ between two symbols $\rho ,
\rho'$ in $S_k$, $k \in {\mathbb Z}$, holds if and only if $\rho -
\rho'$ is a symbol of order $n$ for all integers $n$.
\end{defn}

\bigskip

\begin{defn}
-- The class of pseudo-differential operators is the space $\psi$
modulo addition by an element of $\psi (Z)$, where $Z$ is the
sub-algebra of $S$ with elements equivalent to the zero symbol.
\end{defn}

\bigskip

It is possible to invert the map $\psi$ to obtain for each element
$P$ of $\psi$ a unique symbol $\sigma ( P)$ up to equivalence.
Recall from \S\ref{sec1} that the trace $\tau$ on
$A_{\theta}^{\infty}$ enables one to define the adjoint of operators
acting on $A_{\theta}^{\infty}$ via their extension to ${\cH_0}$. By
direct analogy with \cite{G}, Chapter 1, Theorem, p.~16, one may
deduce the following result.

\bigskip

\begin{prop}
-- For an element $P$ of $\psi$ with symbol $\sigma (P) = \rho =
\rho (\xi)$, the symbol of the adjoint $P^*$ satisfies
\begin{equation*}
\sigma (P^*) \sim \sum_{(\ell_1 , \ell_2) \in ({\mathbb Z}_{\geq
0})^2} \, (1/(\ell_1)! \, (\ell_2)!) \, [\partial_1^{\ell_1} \,
\partial_2^{\ell_2} \, \delta_1^{\ell_1} \, \delta_2^{\ell_2} (\rho
(\xi))^*] \, .
\end{equation*}
If $Q$ is an element of $\psi$ with symbol $\sigma (Q) = \rho' =
\rho' (\xi)$, then the product $PQ$ is also in $\psi$ and has symbol
\begin{equation*}
\sigma (PQ) \sim \sum_{(\ell_1 , \ell_2) \in ({\mathbb Z}_{\geq
0})^2} (1/(\ell_1)! \, \ell_2)!) \, [\partial_1^{\ell_1} \,
\partial_2^{\ell_2} (\rho (\xi)) \, \delta_1^{\ell_1} \,
\delta_2^{\ell_2} (\rho' (\xi))] \, .
\end{equation*}
\end{prop}

\bigskip

Notice that in the above Proposition as throughout the present
paper, given symbols $\{ \rho_j \}_{j=0}^{\infty}$ the relation
$\rho \sim \sum_{j=0}^{\infty} \, \rho_j$ signifies that
for any given $k$ there
exists a positive integer $h$ such that for all $n>h$, the
difference $\rho - \sum_{j=0}^n \, \rho_j$ is in $S_k$.
 The elliptic pseudo-differential operators are those
whose symbols fulfil the criterion which follows.

\bigskip

\begin{defn}
-- Let $n$ be an integer and $\rho$ a symbol of order $n$. Then
$\rho = \rho (\xi)$ is elliptic if it is invertible within the
algebra $C^{\infty} (\R^2 , A_{\theta}^{\infty})$ and if its inverse
satisfies
\begin{equation*}
\Vert \rho (\xi)^{-1} \Vert \leq c(1 + \vert \xi \vert)^{-n}
\end{equation*}
for a constant $c$ depending only on $\rho$ and for $\vert \xi \vert
= (\xi_1^2 + \xi_2^2)^{1/2}$ sufficiently large.

\end{defn}

\bigskip

An example of an elliptic operator is provided by the Laplacian
$\lap = \delta_1^2 + \delta_2^2$ on $A_{\theta}^{\infty}$ introduced
in \S\ref{sec1} which has the corresponding invertible symbol
$\sigma (\lap) = \vert \xi \vert^2$.

\section{Understanding the computations}\label{sec4}

The arguments of this section are kept brief, being direct analogues
of standard ones. Bearing in mind the definition of the zeta
function given in \S\ref{sec2}, we observe that by Cauchy's formula
we have
\begin{equation*}
e^{-t\pert} = (1/2 \pi i) \int_C e^{-t\lambda} (\pert-\lambda {\bf
1})^{-1} \, d\lambda
\end{equation*}
where $\lambda$ is a complex number but not real non-negative, and
$C$ encircles the non-negative real axis in the anti-clockwise
direction without touching it. One then obtains a workable estimate
of $(\pert-\lambda {\bf 1})^{-1}$ by passing to the algebra of
symbols with the important nuance that the scalar $\lambda$ is now
treated as a symbol of order two in all calculations. Using the
definition of a symbol, one can replace the trace in the formula for
the zeta function by an integration in the symbol space (argument
along the diagonal), namely,
\begin{equation*}
\zeta (s) = (1/\Gamma (s)) \int_0^{\infty} \int \tau (\sigma
(e^{-t\pert}) (\xi)) \, t^{s-1} \, d\xi \, dt \, .
\end{equation*}
The function $\Gamma (s)$ has a simple pole at $s=0$ with residue
$1$ so that,
\begin{equation*}
\zeta (0) = {\rm Res}_{s=0} \int_0^{\infty} \int \tau (\sigma
(e^{-t\pert}) (\xi)) \, t^{s-1} \, d\xi \, dt \, .
\end{equation*}
Just as in the arguments employed in the derivation of the
asymptotic formula (see for example \cite{G}),
\begin{equation*}
\int \tau (\sigma (e^{-t\pert}) (\xi)) \, d\xi \sim t^{-1}
\sum_{n=0}^{\infty} \, B_{2n} (\pert) \, t^n \, , \quad t \to 0_+ \,
,
\end{equation*}
one may appeal to the Cauchy formula quoted above. In particular, if
$B_{\lambda}$ denotes (a chosen approximation) to the inverse
operator of $(\lambda {\bf 1} - \pert)$, its symbol has an expansion
of the form
\begin{equation*}
\sigma (B_{\lambda}) = \sigma (B_{\lambda})(\xi) = b_0 (\xi) + b_1
(\xi) + b_2 (\xi) + \ldots
\end{equation*}
where $j$ ranges over the non-negative integers and $b_j (\xi) = b_j
(\xi , \lambda)$ is a symbol of order $-2-j$. As we shall explain at
more length in \S\ref{sec5}, these symbols may be calculated
inductively using the symbol algebra formulae beginning with $b_0
(\xi) = (\lambda - k^2 \vert \xi \vert^2)^{-1}$ which is the
principal (highest homogeneous degree in $\xi$) symbol of $(\lambda
{\bf 1} - \pert)^{-1}$. It turns out that $\zeta (0)$ equals the
coefficient of $\lambda^{-1}$ in $\int \tau (b_2 (\xi)) \, d\xi$. By
a homogeneity argument one has in fact
\begin{equation}
\label{eq1} \zeta (0) = \int \tau (b_2 (\xi)) \, d\xi  \, .
\end{equation}

\section{Computational proof of the main Lemma}\label{sec5}

Following on from the arguments of \S\ref{sec4}, by homogeneity
there is no loss of generality in placing $\lambda = -1$ throughout
the computation of $\zeta (0)$ and multiplying the final answer by
$-$1. The problem is then to derive in the symbol algebra a
recursive solution of the form $\sigma = b_0 (\xi) + b_1 (\xi) + b_2
(\xi) + \ldots$ to the equation
\begin{equation*}
\sigma \cdot (\sigma (\pert-\lambda)) = 1 + 0 (\vert \xi \vert^{-3})
\, .
\end{equation*}
The accuracy to order $-3$ in $\xi$ on the right hand side is in
practice sufficient as we are only interested in evaluating $\sigma$
up to $b_2 (\xi)$. Throughout this section the convention of
summation over repeated indices in the range $i,j = 1,2$ is
observed.

\bigskip

\begin{lem}
-- The operator $\pert$ has symbol $\sigma (\pert) = a_2 (\xi) + a_1
(\xi) + a_0 (\xi)$ where, with summation over repeated indices in
the range $i = 1,2$, one has
\begin{equation*}
a_2 = a_2 (\xi) = k^2 \xi_i \, \xi_i
\end{equation*}
\begin{equation*}
a_1 = a_1 (\xi) = 2 \, \xi_i (k \, \delta_i (k))
\end{equation*}
\begin{equation*}
a_0 = a_0 (\xi) = k \, \delta_i \, \delta_i (k) \, .
\end{equation*}
These expressions are derived by applying the product formula within
the algebra of symbols given in Proposition {\rm \S\ref{sec3}} to
$\sigma_1 (\xi) = \xi_i \, \xi_i$ and $\sigma_2 (\xi) = k$ and then
multiplying on the left by $k$.
\end{lem}

\bigskip

To begin the inductive calculation of the inverse of the symbol of
$\pert-\lambda$, set
\begin{equation}
\label{eq2} b_0 = b_0 (\xi) = (k^2 \, \vert \xi \vert^2 + 1)^{-1}
\end{equation}
and compute to order $-3$ in $\xi$ the product $b_0 \cdot ((a_2 + 1)
+ a_1 + a_0)$. By singling out terms of the appropriate degree $-1$
in $\xi$ and using the Proposition of \S\ref{sec3}, one obtains
\begin{equation}
\label{eq3} b_1 = - (b_0 \, a_1 \, b_0 + \partial_i (b_0) \,
\delta_i (a_2) \, b_0) \, .
\end{equation}
Note that $b_0$ appears on the right in this formula, unlike what is
given in the formula of page 52 of \cite{G} which is valid for
scalar principal symbol only. In a similar fashion, collecting terms
of degree $-2$ in $\xi$ and using (\ref{eq3}) one obtains
\begin{eqnarray}
\label{eq4}
b_2 &= &- (b_0 \, a_0 \, b_0 + b_1 \, a_1 \, b_0 + \partial_i (b_0) \, \delta_i (a_1) \, b_0 \nonumber \\
&+ &\partial_i (b_1) \, \delta_i (a_2) \, b_0 + (1/2) \, \partial_i
\, \partial_j (b_0) \, \delta_i \, \delta_j (a_2) \, b_0) \, .
\end{eqnarray}
A direct computation of these terms (up to the last multiplication
by $b_0$) gives, ignoring the terms which are odd in
$\xi$,\begin{small}
\begin{center}\begin{math}
-b_0k\delta _1^2(k)-b_0k\delta _2^2(k)+\left(\xi_2 ^2+5 \xi_1
^2\right) \left(k^2 b_0^2\right)k\delta _1^2(k)+\left(5 \xi_2
^2+\xi_1 ^2\right) \left(k^2 b_0^2\right)k\delta _2^2(k)+2 \xi_2 ^2
\left(k^2 b_0^2\right)\delta _1(k)\delta _1(k)+6 \xi_1 ^2 \left(k^2
b_0^2\right)\delta _1(k)\delta _1(k)+\xi_2 ^2 \left(k^2
b_0^2\right)\delta _1^2(k)k+\xi_1 ^2 \left(k^2 b_0^2\right)\delta
_1^2(k)k+6 \xi_2 ^2 \left(k^2 b_0^2\right)\delta _2(k)\delta _2(k)+2
\xi_1 ^2 \left(k^2 b_0^2\right)\delta _2(k)\delta _2(k)+\xi_2 ^2
\left(k^2 b_0^2\right)\delta _2^2(k)k+\xi_1 ^2 \left(k^2
b_0^2\right)\delta _2^2(k)k-4 \xi_2 ^2 \xi_1 ^2 \left(k^4
b_0^3\right)k\delta _1^2(k)-4 \xi_1 ^4 \left(k^4 b_0^3\right)k\delta
_1^2(k)-4 \xi_2 ^4 \left(k^4 b_0^3\right)k\delta _2^2(k)-4 \xi_2 ^2
\xi_1 ^2 \left(k^4 b_0^3\right)k\delta _2^2(k)-8 \xi_2 ^2 \xi_1 ^2
\left(k^4 b_0^3\right)\delta _1(k)\delta _1(k)-8 \xi_1 ^4 \left(k^4
b_0^3\right)\delta _1(k)\delta _1(k)-4 \xi_2 ^2 \xi_1 ^2 \left(k^4
b_0^3\right)\delta _1^2(k)k-4 \xi_1 ^4 \left(k^4 b_0^3\right)\delta
_1^2(k)k-8 \xi_2 ^4 \left(k^4 b_0^3\right)\delta _2(k)\delta _2(k)-8
\xi_2 ^2 \xi_1 ^2 \left(k^4 b_0^3\right)\delta _2(k)\delta _2(k)-4
\xi_2 ^4 \left(k^4 b_0^3\right)\delta _2^2(k)k-4 \xi_2 ^2 \xi_1 ^2
\left(k^4 b_0^3\right)\delta _2^2(k)k+2 \xi_2 ^2 b_0k\delta
_1(k)b_0k\delta _1(k)+6 \xi_1 ^2 b_0k\delta _1(k)b_0k\delta _1(k)+2
\xi_2 ^2 b_0k\delta _1(k)b_0\delta _1(k)k+2 \xi_1 ^2 b_0k\delta
_1(k)b_0\delta _1(k)k-4 \xi_2 ^2 \xi_1 ^2 b_0k\delta _1(k)\left(k^2
b_0^2\right)k\delta _1(k)-4 \xi_1 ^4 b_0k\delta _1(k)\left(k^2
b_0^2\right)k\delta _1(k)-4 \xi_2 ^2 \xi_1 ^2 b_0k\delta
_1(k)\left(k^2 b_0^2\right)\delta _1(k)k-4 \xi_1 ^4 b_0k\delta
_1(k)\left(k^2 b_0^2\right)\delta _1(k)k+6 \xi_2 ^2 b_0k\delta
_2(k)b_0k\delta _2(k)+2 \xi_1 ^2 b_0k\delta _2(k)b_0k\delta _2(k)+2
\xi_2 ^2 b_0k\delta _2(k)b_0\delta _2(k)k+2 \xi_1 ^2 b_0k\delta
_2(k)b_0\delta _2(k)k-4 \xi_2 ^4 b_0k\delta _2(k)\left(k^2
b_0^2\right)k\delta _2(k)-4 \xi_2 ^2 \xi_1 ^2 b_0k\delta
_2(k)\left(k^2 b_0^2\right)k\delta _2(k)-4 \xi_2 ^4 b_0k\delta
_2(k)\left(k^2 b_0^2\right)\delta _2(k)k-4 \xi_2 ^2 \xi_1 ^2
b_0k\delta _2(k)\left(k^2 b_0^2\right)\delta _2(k)k-2 \xi_2 ^4
\left(k^2 b_0^2\right)k\delta _1(k)b_0k\delta _1(k)-16 \xi_2 ^2
\xi_1 ^2 \left(k^2 b_0^2\right)k\delta _1(k)b_0k\delta _1(k)-14
\xi_1 ^4 \left(k^2 b_0^2\right)k\delta _1(k)b_0k\delta _1(k)-2 \xi_2
^4 \left(k^2 b_0^2\right)k\delta _1(k)b_0\delta _1(k)k-12 \xi_2 ^2
\xi_1 ^2 \left(k^2 b_0^2\right)k\delta _1(k)b_0\delta _1(k)k-10
\xi_1 ^4 \left(k^2 b_0^2\right)k\delta _1(k)b_0\delta _1(k)k+4 \xi_2
^4 \xi_1 ^2 \left(k^2 b_0^2\right)k\delta _1(k)\left(k^2
b_0^2\right)k\delta _1(k)+8 \xi_2 ^2 \xi_1 ^4 \left(k^2
b_0^2\right)k\delta _1(k)\left(k^2 b_0^2\right)k\delta _1(k)+4 \xi_1
^6 \left(k^2 b_0^2\right)k\delta _1(k)\left(k^2 b_0^2\right)k\delta
_1(k)+4 \xi_2 ^4 \xi_1 ^2 \left(k^2 b_0^2\right)k\delta
_1(k)\left(k^2 b_0^2\right)\delta _1(k)k+8 \xi_2 ^2 \xi_1 ^4
\left(k^2 b_0^2\right)k\delta _1(k)\left(k^2 b_0^2\right)\delta
_1(k)k+4 \xi_1 ^6 \left(k^2 b_0^2\right)k\delta _1(k)\left(k^2
b_0^2\right)\delta _1(k)k-14 \xi_2 ^4 \left(k^2 b_0^2\right)k\delta
_2(k)b_0k\delta _2(k)-16 \xi_2 ^2 \xi_1 ^2 \left(k^2
b_0^2\right)k\delta _2(k)b_0k\delta _2(k)-2 \xi_1 ^4 \left(k^2
b_0^2\right)k\delta _2(k)b_0k\delta _2(k)-10 \xi_2 ^4 \left(k^2
b_0^2\right)k\delta _2(k)b_0\delta _2(k)k-12 \xi_2 ^2 \xi_1 ^2
\left(k^2 b_0^2\right)k\delta _2(k)b_0\delta _2(k)k-2 \xi_1 ^4
\left(k^2 b_0^2\right)k\delta _2(k)b_0\delta _2(k)k+4 \xi_2 ^6
\left(k^2 b_0^2\right)k\delta _2(k)\left(k^2 b_0^2\right)k\delta
_2(k)+8 \xi_2 ^4 \xi_1 ^2 \left(k^2 b_0^2\right)k\delta
_2(k)\left(k^2 b_0^2\right)k\delta _2(k)+4 \xi_2 ^2 \xi_1 ^4
\left(k^2 b_0^2\right)k\delta _2(k)\left(k^2 b_0^2\right)k\delta
_2(k)+4 \xi_2 ^6 \left(k^2 b_0^2\right)k\delta _2(k)\left(k^2
b_0^2\right)\delta _2(k)k+8 \xi_2 ^4 \xi_1 ^2 \left(k^2
b_0^2\right)k\delta _2(k)\left(k^2 b_0^2\right)\delta _2(k)k+4 \xi_2
^2 \xi_1 ^4 \left(k^2 b_0^2\right)k\delta _2(k)\left(k^2
b_0^2\right)\delta _2(k)k-2 \xi_2 ^4 \left(k^2 b_0^2\right)\delta
_1(k)kb_0k\delta _1(k)-12 \xi_2 ^2 \xi_1 ^2 \left(k^2
b_0^2\right)\delta _1(k)kb_0k\delta _1(k)-10 \xi_1 ^4 \left(k^2
b_0^2\right)\delta _1(k)kb_0k\delta _1(k)-2 \xi_2 ^4 \left(k^2
b_0^2\right)\delta _1(k)kb_0\delta _1(k)k-8 \xi_2 ^2 \xi_1 ^2
\left(k^2 b_0^2\right)\delta _1(k)kb_0\delta _1(k)k-6 \xi_1 ^4
\left(k^2 b_0^2\right)\delta _1(k)kb_0\delta _1(k)k+4 \xi_2 ^4 \xi_1
^2 \left(k^2 b_0^2\right)\delta _1(k)k\left(k^2 b_0^2\right)k\delta
_1(k)+8 \xi_2 ^2 \xi_1 ^4 \left(k^2 b_0^2\right)\delta
_1(k)k\left(k^2 b_0^2\right)k\delta _1(k)+4 \xi_1 ^6 \left(k^2
b_0^2\right)\delta _1(k)k\left(k^2 b_0^2\right)k\delta _1(k)+4 \xi_2
^4 \xi_1 ^2 \left(k^2 b_0^2\right)\delta _1(k)k\left(k^2
b_0^2\right)\delta _1(k)k+8 \xi_2 ^2 \xi_1 ^4 \left(k^2
b_0^2\right)\delta _1(k)k\left(k^2 b_0^2\right)\delta _1(k)k+4 \xi_1
^6 \left(k^2 b_0^2\right)\delta _1(k)k\left(k^2 b_0^2\right)\delta
_1(k)k-10 \xi_2 ^4 \left(k^2 b_0^2\right)\delta _2(k)kb_0k\delta
_2(k)-12 \xi_2 ^2 \xi_1 ^2 \left(k^2 b_0^2\right)\delta
_2(k)kb_0k\delta _2(k)-2 \xi_1 ^4 \left(k^2 b_0^2\right)\delta
_2(k)kb_0k\delta _2(k)-6 \xi_2 ^4 \left(k^2 b_0^2\right)\delta
_2(k)kb_0\delta _2(k)k-8 \xi_2 ^2 \xi_1 ^2 \left(k^2
b_0^2\right)\delta _2(k)kb_0\delta _2(k)k-2 \xi_1 ^4 \left(k^2
b_0^2\right)\delta _2(k)kb_0\delta _2(k)k+4 \xi_2 ^6 \left(k^2
b_0^2\right)\delta _2(k)k\left(k^2 b_0^2\right)k\delta _2(k)+8 \xi_2
^4 \xi_1 ^2 \left(k^2 b_0^2\right)\delta _2(k)k\left(k^2
b_0^2\right)k\delta _2(k)+4 \xi_2 ^2 \xi_1 ^4 \left(k^2
b_0^2\right)\delta _2(k)k\left(k^2 b_0^2\right)k\delta _2(k)+4 \xi_2
^6 \left(k^2 b_0^2\right)\delta _2(k)k\left(k^2 b_0^2\right)\delta
_2(k)k+8 \xi_2 ^4 \xi_1 ^2 \left(k^2 b_0^2\right)\delta
_2(k)k\left(k^2 b_0^2\right)\delta _2(k)k+4 \xi_2 ^2 \xi_1 ^4
\left(k^2 b_0^2\right)\delta _2(k)k\left(k^2 b_0^2\right)\delta
_2(k)k+8 \xi_2 ^4 \xi_1 ^2 \left(k^4 b_0^3\right)k\delta
_1(k)b_0k\delta _1(k)+16 \xi_2 ^2 \xi_1 ^4 \left(k^4
b_0^3\right)k\delta _1(k)b_0k\delta _1(k)+8 \xi_1 ^6 \left(k^4
b_0^3\right)k\delta _1(k)b_0k\delta _1(k)+8 \xi_2 ^4 \xi_1 ^2
\left(k^4 b_0^3\right)k\delta _1(k)b_0\delta _1(k)k+16 \xi_2 ^2
\xi_1 ^4 \left(k^4 b_0^3\right)k\delta _1(k)b_0\delta _1(k)k+8 \xi_1
^6 \left(k^4 b_0^3\right)k\delta _1(k)b_0\delta _1(k)k+8 \xi_2 ^6
\left(k^4 b_0^3\right)k\delta _2(k)b_0k\delta _2(k)+16 \xi_2 ^4
\xi_1 ^2 \left(k^4 b_0^3\right)k\delta _2(k)b_0k\delta _2(k)+8 \xi_2
^2 \xi_1 ^4 \left(k^4 b_0^3\right)k\delta _2(k)b_0k\delta _2(k)+8
\xi_2 ^6 \left(k^4 b_0^3\right)k\delta _2(k)b_0\delta _2(k)k+16
\xi_2 ^4 \xi_1 ^2 \left(k^4 b_0^3\right)k\delta _2(k)b_0\delta
_2(k)k+8 \xi_2 ^2 \xi_1 ^4 \left(k^4 b_0^3\right)k\delta
_2(k)b_0\delta _2(k)k+8 \xi_2 ^4 \xi_1 ^2 \left(k^4
b_0^3\right)\delta _1(k)kb_0k\delta _1(k)+16 \xi_2 ^2 \xi_1 ^4
\left(k^4 b_0^3\right)\delta _1(k)kb_0k\delta _1(k)+8 \xi_1 ^6
\left(k^4 b_0^3\right)\delta _1(k)kb_0k\delta _1(k)+8 \xi_2 ^4 \xi_1
^2 \left(k^4 b_0^3\right)\delta _1(k)kb_0\delta _1(k)k+16 \xi_2 ^2
\xi_1 ^4 \left(k^4 b_0^3\right)\delta _1(k)kb_0\delta _1(k)k+8 \xi_1
^6 \left(k^4 b_0^3\right)\delta _1(k)kb_0\delta _1(k)k+8 \xi_2 ^6
\left(k^4 b_0^3\right)\delta _2(k)kb_0k\delta _2(k)+16 \xi_2 ^4
\xi_1 ^2 \left(k^4 b_0^3\right)\delta _2(k)kb_0k\delta _2(k)+8 \xi_2
^2 \xi_1 ^4 \left(k^4 b_0^3\right)\delta _2(k)kb_0k\delta _2(k)+8
\xi_2 ^6 \left(k^4 b_0^3\right)\delta _2(k)kb_0\delta _2(k)k+16
\xi_2 ^4 \xi_1 ^2 \left(k^4 b_0^3\right)\delta _2(k)kb_0\delta
_2(k)k+8 \xi_2 ^2 \xi_1 ^4 \left(k^4 b_0^3\right)\delta
_2(k)kb_0\delta _2(k)k
\end{math}\end{center}
\end{small}
It is extremely useful during the computation to exploit the fact
that in the target formula for $\zeta (0)$ given in (\ref{eq1}),
\S\ref{sec5}, one invokes the trace, so that members of the factors
of the individual summands may be permuted cyclically without loss
of generality for the answer. In our case this means that instead of
multiplying the above sum by $b_0$ on the right we can just multiply
it on the left, thus simply adding one to the exponent of the first
occurrence of $b_0$. By formula (\ref{eq1}), one then has to sum the
integrals of each of these terms over the whole $\xi$-plane. After
multiplying by $b_0$ on the left, passing in polar coordinates and
integrating the angular variable one gets, up to an overall factor
of $2\pi$,\begin{small}
\begin{center}\begin{math}
-b_0^2k\delta _1^2(k)-b_0^2k\delta _2^2(k)+3 r^2 \left(k^2
b_0^3\right)k\delta _1^2(k)+3 r^2 \left(k^2 b_0^3\right)k\delta
_2^2(k)+4 r^2 \left(k^2 b_0^3\right)\delta _1(k)\delta _1(k)+r^2
\left(k^2 b_0^3\right)\delta _1^2(k)k+4 r^2 \left(k^2
b_0^3\right)\delta _2(k)\delta _2(k)+r^2 \left(k^2
b_0^3\right)\delta _2^2(k)k-2 r^4 \left(k^4 b_0^4\right)k\delta
_1^2(k)-2 r^4 \left(k^4 b_0^4\right)k\delta _2^2(k)-4 r^4 \left(k^4
b_0^4\right)\delta _1(k)\delta _1(k)-2 r^4 \left(k^4
b_0^4\right)\delta _1^2(k)k-4 r^4 \left(k^4 b_0^4\right)\delta
_2(k)\delta _2(k)-2 r^4 \left(k^4 b_0^4\right)\delta _2^2(k)k+4 r^2
b_0^2k\delta _1(k)b_0k\delta _1(k)+2 r^2 b_0^2k\delta _1(k)b_0\delta
_1(k)k-2 r^4 b_0^2k\delta _1(k)\left(k^2 b_0^2\right)k\delta _1(k)-2
r^4 b_0^2k\delta _1(k)\left(k^2 b_0^2\right)\delta _1(k)k+4 r^2
b_0^2k\delta _2(k)b_0k\delta _2(k)+2 r^2 b_0^2k\delta _2(k)b_0\delta
_2(k)k-2 r^4 b_0^2k\delta _2(k)\left(k^2 b_0^2\right)k\delta _2(k)-2
r^4 b_0^2k\delta _2(k)\left(k^2 b_0^2\right)\delta _2(k)k-8 r^4
\left(k^2 b_0^3\right)k\delta _1(k)b_0k\delta _1(k)-6 r^4 \left(k^2
b_0^3\right)k\delta _1(k)b_0\delta _1(k)k+2 r^6 \left(k^2
b_0^3\right)k\delta _1(k)\left(k^2 b_0^2\right)k\delta _1(k)+2 r^6
\left(k^2 b_0^3\right)k\delta _1(k)\left(k^2 b_0^2\right)\delta
_1(k)k-8 r^4 \left(k^2 b_0^3\right)k\delta _2(k)b_0k\delta _2(k)-6
r^4 \left(k^2 b_0^3\right)k\delta _2(k)b_0\delta _2(k)k+2 r^6
\left(k^2 b_0^3\right)k\delta _2(k)\left(k^2 b_0^2\right)k\delta
_2(k)+2 r^6 \left(k^2 b_0^3\right)k\delta _2(k)\left(k^2
b_0^2\right)\delta _2(k)k-6 r^4 \left(k^2 b_0^3\right)\delta
_1(k)kb_0k\delta _1(k)-4 r^4 \left(k^2 b_0^3\right)\delta
_1(k)kb_0\delta _1(k)k+2 r^6 \left(k^2 b_0^3\right)\delta
_1(k)k\left(k^2 b_0^2\right)k\delta _1(k)+2 r^6 \left(k^2
b_0^3\right)\delta _1(k)k\left(k^2 b_0^2\right)\delta _1(k)k-6 r^4
\left(k^2 b_0^3\right)\delta _2(k)kb_0k\delta _2(k)-4 r^4 \left(k^2
b_0^3\right)\delta _2(k)kb_0\delta _2(k)k+2 r^6 \left(k^2
b_0^3\right)\delta _2(k)k\left(k^2 b_0^2\right)k\delta _2(k)+2 r^6
\left(k^2 b_0^3\right)\delta _2(k)k\left(k^2 b_0^2\right)\delta
_2(k)k+4 r^6 \left(k^4 b_0^4\right)k\delta _1(k)b_0k\delta _1(k)+4
r^6 \left(k^4 b_0^4\right)k\delta _1(k)b_0\delta _1(k)k+4 r^6
\left(k^4 b_0^4\right)k\delta _2(k)b_0k\delta _2(k)+4 r^6 \left(k^4
b_0^4\right)k\delta _2(k)b_0\delta _2(k)k+4 r^6 \left(k^4
b_0^4\right)\delta _1(k)kb_0k\delta _1(k)+4 r^6 \left(k^4
b_0^4\right)\delta _1(k)kb_0\delta _1(k)k+4 r^6 \left(k^4
b_0^4\right)\delta _2(k)kb_0k\delta _2(k)+4 r^6 \left(k^4
b_0^4\right)\delta _2(k)kb_0\delta _2(k)k
\end{math}\end{center}\end{small}
where $b_0=(r^2k^2-\lambda)^{-1}$ and where the integration is in
$rdr$ and from $0$ to $\infty$.

\subsection{Terms with  all $b_0$ on the left}
Using the trace property these terms give the following:
\begin{center}\begin{math}
-b_0^2 k\delta _1^2(k)-b_0^2 k\delta _2^2(k)+r^2 (3 k^2 b_0^3
k\delta _1^2(k)-2 k^4 r^2 b_0^4 k\delta _1^2(k)+3 k^2 b_0^3 k\delta
_2^2(k)-2 k^4 r^2 b_0^4 k\delta _2^2(k)+4 k^2 b_0^3 \delta
_1(k)\delta _1(k)-4 k^4 r^2 b_0^4 \delta _1(k)\delta _1(k)+k^2 b_0^3
\delta _1^2(k)k-2 k^4 r^2 b_0^4 \delta _1^2(k)k+4 k^2 b_0^3 \delta
_2(k)\delta _2(k)-4 k^4 r^2 b_0^4 \delta _2(k)\delta _2(k)+k^2 b_0^3
\delta _2^2(k)k-2 k^4 r^2 b_0^4 \delta _2^2(k)k)
\end{math}\end{center}
which gives the same result as
\begin{center}\begin{math}
\left(4 k^2 r^2 b_0^3-4 k^4 r^4 b_0^4\right) \left(\delta
_1(k){}^2+\delta _2(k){}^2\right)+\left(-k b_0^2+4 k^3 r^2 b_0^3-4
k^5 r^4 b_0^4\right) \left(\delta _1^2(k)+\delta _2^2(k)\right)
\end{math}\end{center}
and the coefficient of $1/\lambda$ in the integral $\int \bullet
\,\,rdr$ gives, up to an overall coefficient of $2\pi$,
\begin{equation}\label{res1}
    -\frac 13 k^{-2}\left(\delta _1(k){}^2+\delta _2(k){}^2\right)+\frac 16 k^{-1}\left(\delta _1^2(k)+\delta _2^2(k)\right)
\end{equation}
which corresponds to the first two terms of the formula for
$h(\theta,k)$ in the statement of the Lemma in the form
\eqref{former}.

\subsection{Terms with $b_0^2$ in the middle}\label{termT}
They are the following
\begin{center}\begin{math}
-2 r^4 b_0^2k\delta _1(k)\left(k^2 b_0^2\right)k\delta _1(k)-2 r^4
b_0^2k\delta _1(k)\left(k^2 b_0^2\right)\delta _1(k)k-2 r^4
b_0^2k\delta _2(k)\left(k^2 b_0^2\right)k\delta _2(k)-2 r^4
b_0^2k\delta _2(k)\left(k^2 b_0^2\right)\delta _2(k)k+2 r^6
\left(k^2 b_0^3\right)k\delta _1(k)\left(k^2 b_0^2\right)k\delta
_1(k)+2 r^6 \left(k^2 b_0^3\right)k\delta _1(k)\left(k^2
b_0^2\right)\delta _1(k)k+2 r^6 \left(k^2 b_0^3\right)k\delta
_2(k)\left(k^2 b_0^2\right)k\delta _2(k)+2 r^6 \left(k^2
b_0^3\right)k\delta _2(k)\left(k^2 b_0^2\right)\delta _2(k)k+2 r^6
\left(k^2 b_0^3\right)\delta _1(k)k\left(k^2 b_0^2\right)k\delta
_1(k)+2 r^6 \left(k^2 b_0^3\right)\delta _1(k)k\left(k^2
b_0^2\right)\delta _1(k)k+2 r^6 \left(k^2 b_0^3\right)\delta
_2(k)k\left(k^2 b_0^2\right)k\delta _2(k)+2 r^6 \left(k^2
b_0^3\right)\delta _2(k)k\left(k^2 b_0^2\right)\delta _2(k)k
\end{math}\end{center}

One has $$
\partial_r(b_0)=-2k^2 r b_0^2
$$
and one can use integration by parts in $r$ to transform terms such
as
$$
\int_0^\infty r^6 (k^2 b_0^3)\delta _1(k)k(k^2 b_0^2)\delta _1(k)k
rdr
$$
and replace them by
$$
\frac 12\int_0^\infty \partial_r\left(r^6 (k^2 b_0^3)\right)\delta
_1(k)k b_0\delta _1(k)k dr
$$
where now $b_0$ only appears at the first power on the second
occurrence.

After performing this operation, and combining with those which have
$b_0$ in the middle, namely
\begin{center}\begin{math}
4 r^2 b_0^2k\delta _1(k)b_0k\delta _1(k)+2 r^2 b_0^2k\delta
_1(k)b_0\delta _1(k)k+4 r^2 b_0^2k\delta _2(k)b_0k\delta _2(k)+2 r^2
b_0^2k\delta _2(k)b_0\delta _2(k)k-8 r^4 \left(k^2
b_0^3\right)k\delta _1(k)b_0k\delta _1(k)-6 r^4 \left(k^2
b_0^3\right)k\delta _1(k)b_0\delta _1(k)k-8 r^4 \left(k^2
b_0^3\right)k\delta _2(k)b_0k\delta _2(k)-6 r^4 \left(k^2
b_0^3\right)k\delta _2(k)b_0\delta _2(k)k-6 r^4 \left(k^2
b_0^3\right)\delta _1(k)kb_0k\delta _1(k)-4 r^4 \left(k^2
b_0^3\right)\delta _1(k)kb_0\delta _1(k)k-6 r^4 \left(k^2
b_0^3\right)\delta _2(k)kb_0k\delta _2(k)-4 r^4 \left(k^2
b_0^3\right)\delta _2(k)kb_0\delta _2(k)k+4 r^6 \left(k^4
b_0^4\right)k\delta _1(k)b_0k\delta _1(k)+4 r^6 \left(k^4
b_0^4\right)k\delta _1(k)b_0\delta _1(k)k+4 r^6 \left(k^4
b_0^4\right)k\delta _2(k)b_0k\delta _2(k)+4 r^6 \left(k^4
b_0^4\right)k\delta _2(k)b_0\delta _2(k)k+4 r^6 \left(k^4
b_0^4\right)\delta _1(k)kb_0k\delta _1(k)+4 r^6 \left(k^4
b_0^4\right)\delta _1(k)kb_0\delta _1(k)k+4 r^6 \left(k^4
b_0^4\right)\delta _2(k)kb_0k\delta _2(k)+4 r^6 \left(k^4
b_0^4\right)\delta _2(k)kb_0\delta _2(k)k
\end{math}\end{center}

one obtains the following terms
\begin{center}\begin{math}
T=-2 \left(r^2 b_0^2\right)\left(k\delta _1(k)b_0\delta
_1(k)k+k\delta _2(k)b_0\delta _2(k)k\right)+2 \left(k^2 r^4
b_0^3\right)(k\delta _1(k)b_0k\delta _1(k)+2 k\delta _1(k)b_0\delta
_1(k)k+k\delta _2(k)b_0k\delta _2(k)+2 k\delta _2(k)b_0\delta
_2(k)k+\delta _1(k)kb_0\delta _1(k)k+\delta _2(k)kb_0\delta
_2(k)k)-2 \left(k^4 r^6 b_0^4\right)(k\delta _1(k)b_0k\delta
_1(k)+k\delta _1(k)b_0\delta _1(k)k+k\delta _2(k)b_0k\delta
_2(k)+k\delta _2(k)b_0\delta _2(k)k+\delta _1(k)kb_0k\delta
_1(k)+\delta _1(k)kb_0\delta _1(k)k+\delta _2(k)kb_0k\delta
_2(k)+\delta _2(k)kb_0\delta _2(k)k)
\end{math}\end{center}
which all have $b_0$ in the middle.

\subsection{Terms with $b_0$ in the middle}

 Since we are in the non-commutative case, when in particular $k$ and $\delta_i ( k)$, $i= 1,2$, do not commute, the
computation of such terms requires us to permute $k$ with elements of
$A_{\theta}^{\infty}$. This is achieved using the following lemma.

\medskip

\begin{lem}\label{permlem}
-- For every element $\rho$ of $A_{\theta}^{\infty}$ and every
non-negative integer $m$ one has,
\begin{equation}\label{move}
\int_0^{\infty} \frac{k^{2m+2} \, u^m }{ (k^2 \, u + 1)^{m+1}}
\,\rho\, \frac{1}{(k^2 \, u+1)} \, du = {\mathcal D}_m (\rho) \,,
\end{equation}
where the modified logarithm function is ${\mathcal D}_m={\mathcal
L}_m (\modu)$, $\modu$ is the operator introduced in \S\ref{sec1},
and
\begin{equation}\label{move1}
{\mathcal L}_m (u)= \int_0^{\infty} \frac{x^m }{
(x+1)^{m+1}}\frac{1}{(xu + 1)} \, dx \, .
\end{equation}
\end{lem}

\proof On effecting the change of variables $u=e^s$ one obtains,
with $k=e^{f/2}$,
\begin{eqnarray}
&&\int_0^{\infty} \frac{k^{2m+2} \, u^m }{ (k^2 \, u + 1)^{m+1}} \,\rho\, \frac{1}{(k^2 \, u+1)} \, du \nonumber \\
&= &\int_{-\infty}^{\infty} \frac{e^{(m+1)f+ms} }{ (e^{(s+f)} + 1)^{m+1}} \,\rho\, \frac{e^s}{(e^{(s+f)} + 1)} \, ds \nonumber \\
&= &\int_{-\infty}^{\infty} \frac{e^{(m+1/2)(s+f)} }{ (e^{(s+f)} + 1)^{m+1}} \,\modu^{-1/2} (\rho)\, \frac{e^{(s+f) / 2}}{(e^{(s+f)} + 1)} \, ds \nonumber \\
&= & \int_{-\infty}^{\infty} \frac{e^{(m+1/2)(s+f)} }{ (e^{(s+f)} + 1)^{m+1}} \,\modu^{-1/2} (\rho)\, \int_{-\infty}^{\infty} \frac{e^{it(s+f)} }{e^{\pi t} + e^{-\pi t}} \, dt \, ds \nonumber \\
&= & \int_{-\infty}^{\infty} \frac{e^{(m+1/2)(s+f)} }{ (e^{(s+f)} + 1)^{m+1}} \, \int_{-\infty}^{\infty} \frac{e^{it(s+f)} }{e^{\pi t} + e^{-\pi t}} \, \modu^{-1/2+it} (\rho) \,  dt \, ds \nonumber \\
&= & \int_{-\infty}^{\infty} \frac{e^{itf} }{e^{\pi t} + e^{-\pi
t}}\left( \int_{-\infty}^{\infty} \frac{e^{(m+1/2)(s+f)} e^{its}}{
(e^{(s+f)} + 1)^{m+1}}  \, ds\right)\,  \modu^{-1/2+it} (\rho) \,
dt \nonumber
\end{eqnarray}
Now one has
$$
\int_{-\infty}^{\infty} \frac{e^{(m+1/2)(s+f)} e^{its}}{ (e^{(s+f)}
+ 1)^{m+1}}  \, ds=e^{-itf}F_m(t)
$$
where $F_m$ is the Fourier transform of the function
$$h_m(s)=\frac{e^{(m+1/2)s} }{ (e^{s} + 1)^{m+1}} $$ We thus get
\begin{eqnarray}
&&\int_0^{\infty} \frac{k^{2m+2} \, u^m }{ (k^2 \, u + 1)^{m+1}} \,\rho\, \frac{1}{(k^2 \, u+1)} \, du \nonumber \\
&= &\int_{-\infty}^{\infty} \frac{F_m(t) }{e^{\pi t} + e^{-\pi t}}\,
\modu^{-1/2+it} (\rho) \,  dt \nonumber  \, . \nonumber
\end{eqnarray}
Moreover one has
\begin{eqnarray}
&&\int_{-\infty}^{\infty} \frac{F_m(t) }{e^{\pi t} + e^{-\pi t}}\,  u^{-it} \,  dt \nonumber \\
&= &\int_0^{\infty} \frac{x^{(m+1/2)} }{ (x+1)^{m+1}}
\frac{x^{-1/2}u^{1/2}}{x^{-1}u +1} \, \frac{dx }{x} =
u^{-1/2}\,{\mathcal L}_m (u^{-1}) \, . \nonumber
\end{eqnarray}
Replacing $u$ by $\modu^{-1}$ one gets the required equality
\eqref{move}. One can check the normalization taking $u=1$. Then the
integral ${\mathcal L}_m$ equals $I_m$
 \begin{equation*}
I_m = \int_0^{\infty} v^m / (v+1)^{m+2} \, dv = 1/(m+1) \, ,
\end{equation*}
for every positive integer $m$. On the other hand, by inspection one
sees that ${\mathcal L}_m$ is of the form
\begin{equation*}
{\mathcal L}_m(u) = c_m  \left(\log (u) - P (u) \right)/ (u -
1)^{m+1} \, ,
\end{equation*}
where $P$ is a polynomial of degree at most $m$. In the neighborhood
of $u = 1$ one has
\begin{equation*}
\log (u) = \sum_{j=1}^{\infty} \frac{(-1)^{j+1} }{ j} (u - 1)^j \, ,
\end{equation*}
and from   its value at $u = 1$, where ${\mathcal L}_m$ is
non-singular, one sees from this last expression that ${\mathcal
L}_m$ is the modified logarithm ${\mathcal D}_m$ introduced in
\S\ref{sec2} where,
\begin{equation*}
{\mathcal D}_m(\modu) = ((-1)^m / (\modu - 1)^{m+1}) \left\{ \log
(\modu) - \sum_{j=1}^m  \frac{(-1)^{j+1} }{ j} (\modu - 1)^j
\right\} \, .
\end{equation*}
This completes the proof of the Lemma. \endproof

\bigskip

We now split the sum $T$ of \S \ref{termT} as a sum of three terms
$T=T_1+T_2+T_3$ and compute the coefficient of $1/\lambda$ in the
integral with respect to $rdr$ in each of them using the above
Lemma.

\subsubsection{Terms involving ${\mathcal D}_1$}

These terms come from
$$
T_1=-2 \left(r^2 b_0^2\right)\left(k\delta _1(k)b_0\delta
_1(k)k+k\delta _2(k)b_0\delta _2(k)k\right)
$$
With $u=r^2$ the integrand is $rdr=\frac 12 du$ and thus (up to the
overall factor of $2\pi$) these terms give, by setting $\lambda=-1$
and changing the overall sign,
\begin{equation}\label{d1term}
\tau(k^{-2}\int_0^\infty \frac{k^{4} \, u }{ (k^2 \, u + 1)^{2}}
\,\delta_i (k)\, \frac{1}{(k^2 \, u+1)} \, du\,\delta_i (k)) =
\tau({\mathcal D}_1  (\delta_i (k)) \, \delta_i (k) \, k^{-2})
\end{equation}

\subsubsection{Terms involving ${\mathcal D}_2$}

These terms come from
\begin{center}\begin{math}
T_2=2 \left(k^2 r^4 b_0^3\right)(k\delta _1(k)b_0k\delta _1(k)+2
k\delta _1(k)b_0\delta _1(k)k+k\delta _2(k)b_0k\delta _2(k)+2
k\delta _2(k)b_0\delta _2(k)k+\delta _1(k)kb_0\delta _1(k)k+\delta
_2(k)kb_0\delta _2(k)k)
\end{math}\end{center}
Since $k$ commutes with $b_0$ and one works under the trace, they
give the same  as
$$
4 \left(k^2 r^4 b_0^3\right)(k\delta _j(k)k b_0\delta
_j(k)+k^2\delta _j(k)b_0\delta _j(k))
$$
One has $k\delta _j(k)k=k^2\modu^{1/2}(\delta _j(k))$. Thus after
setting $\lambda=-1$ and changing the overall sign,  one gets
$$
-2\tau(k^{-2}\int_0^\infty \frac{k^{6} \, u^2 }{ (k^2 \, u + 1)^{3}}
\,(\modu^{1/2}(\delta _i(k))+\delta_i (k))\, \frac{1}{(k^2 \, u+1)}
\, du\,\delta_i (k))
$$
\begin{equation}\label{d2term}
=-2\tau(({\mathcal D}_2 \, \modu^{1/2})(\delta_i (k)) \, \delta_i
(k) \, k^{-2})-2\tau({\mathcal D}_2 (\delta_i (k)) \, \delta_i (k)
\, k^{-2})
\end{equation}

\subsubsection{Terms involving ${\mathcal D}_3$}

These terms come from
\begin{center}\begin{math}
T_3=-2 \left(k^4 r^6 b_0^4\right)(k\delta _1(k)b_0k\delta
_1(k)+k\delta _1(k)b_0\delta _1(k)k+k\delta _2(k)b_0k\delta
_2(k)+k\delta _2(k)b_0\delta _2(k)k+\delta _1(k)kb_0k\delta
_1(k)+\delta _1(k)kb_0\delta _1(k)k+\delta _2(k)kb_0k\delta
_2(k)+\delta _2(k)kb_0\delta _2(k)k)
\end{math}\end{center}
Since $k$ commutes with $b_0$ and one works under the trace, they
give the same  as
$$
-2 \left(k^4 r^6 b_0^4\right)\left(k^2\delta _j(k)b_0\delta
_j(k)+2k\delta _j(k)k b_0\delta _j(k)+\delta _j(k)k^2 b_0\delta
_j(k)\right)
$$
One has $k\delta _j(k)k=k^2\modu^{1/2}(\delta _j(k))$ and $\delta
_j(k)k^2=k^2\modu(\delta _j(k))$. Thus after setting $\lambda=-1$
and changing the overall sign,  one gets
\begin{equation}\label{d3term}
\tau\left(\left({\mathcal D}_3(\delta_i (k)) \, \delta_i (k)  +2
({\mathcal D}_3 \, \modu^{1/2})(\delta_i (k)) \, \delta_i (k) +
({\mathcal D}_3 \, \modu)(\delta_i (k)) \, \delta_i (k)\right)
k^{-2}\right)
\end{equation}

\end{document}